\documentclass[12pt]{article}
\usepackage{ct}

\usepackage{amsmath, amsthm, amssymb, amsbsy, amsfonts, bbm, bigints, color, comment, float, graphics, graphicx, hyperref,import, mathrsfs, mathtools,picture, stmaryrd, subfig, svg, xcolor}

\usepackage[percent]{overpic}

\newtheorem{thm}{Theorem}[section]

\newtheorem{lem}[thm]{Lemma}
\newtheorem{prop}[thm]{Proposition}

\theoremstyle{definition}
\newtheorem{defn}[thm]{Definition}
\newtheorem{rem}[thm]{Remark}

\theoremstyle{remark}

\newcommand{\R}{\mathbb{R}}
\newcommand{\C}{\mathbb{C}}
\newcommand{\Z}{\mathbb{Z}}
\newcommand{\N}{\mathbb{N}}

\renewcommand{\P}{\mathbb{P}}

\newcommand{\BB}{\textbf{B}}

\newcommand{\ee}{\textbf{e}}
\renewcommand{\aa}{\textbf{a}}
\newcommand{\kk}{\textbf{k}}
\newcommand{\nn}{\textbf{n}}
\newcommand{\rr}{\textbf{r}}

\newcommand{\vv}{\textbf{v}}
\newcommand{\xx}{\textbf{x}}
\newcommand{\yy}{\textbf{y}}
\newcommand{\ttheta}{\boldsymbol{\theta}}
\newcommand{\aalpha}{\boldsymbol{\alpha}}
\newcommand{\un}{\boldsymbol{1}}

\specs{n}{vol (issue)}{year}{}

\dateline{February 3, 2023}{}{}



\title{The asymptotic number of lattice zonotopes in a hypercube.}

\author[1,2]{Buffière, Théophile}

\affil[1]{LIPN, Université Sorbonne Paris Nord (P13), 93430 Villetaneuse, France \newline
\email{buffiere@math.univ-paris13.fr}%
}

\affil[2]{LAGA, Université Sorbonne Paris Nord (P13), 93430 Villetaneuse, France\newline
\email{buffiere@math.univ-paris13.fr}%
}

\begin{document}

\maketitle


\begin{abstract}
We provide a sharp estimate for the asymptotic number of lattice zonotopes, inscribed in $[0,n ]^d$ when $n$ tends to infinity. Our estimate refines the logarithmic equivalent established by Barany, Bureaux, and Lund when the sum of the generators of the zonotope is prescribed. As we shall see, the exponential part of our estimate is composed of a polynomial of degree $d$ in $n^{1/(d+1)}$, and involves Riemann's zeta function and its non-trivial zeros. 
We also analyze some combinatorial properties of lattice zonotopes. In particular, we provide the first moment of the polyhedral graph asymptotic diameter when $n$ goes to infinity.\medskip
\end{abstract}




\section{INTRODUCTION}\label{introduction}

The idea of enumerating geometric objects such as convex integral polytopes is a long-standing question: in the late '70s, Arnold bounded the number of convex lattice polygons with area $A$ up to the affine transformations of $\Z^2$  \cite{Arnold:polygons}, showing that it behaves roughly like $A^{1/3}$.  In 1992, B\'ar\'any and Pach narrowed those bounds \cite{Barany:numberpolygons}, before B\'ar\'any and Vershik  extended the result to higher dimensions \cite{Barany:numberpolytope}. These bounds had, in fact, partially been improved by Konyagin and Sevast'yanov as of 1984 \cite{Konyakin:bound_polyhedron}.
 Vershik later simply raised the question of the number of lattice polytopes in a box \cite{Vershik:limit}; e.g. in 2 dimensions, the number of convex lattice polygons inscribed in the square $[-n, n]^2$. \medskip

A natural intermediate stage is to enumerate the number $p(n)$ of convex polygonal lines from $(0,0)$ to $(n,n)$ on the lattice $\Z^2$ (referred later as lattice chains). In 1994, an asymptotic equivalent of its logarithm is independently found by B\'ar\'any \cite{Barany:limitshape}, Vershik \cite{Vershik:limit}, and Sinai \cite{Sinai:polygonallines}, who showed that:

\begin{align}\label{kappa}
    p(n) = \exp \left(3 \kappa^{1/3} n^{2/3} (1 + \text{o}(1))\right), \:\:\:\:\:\: \text{where} \:\: \kappa = \frac{\zeta(3)}{\zeta(2)} .
\end{align}

More recently, Bodini, Duchon, Jacquot, and Mutafchiev \cite{Bodini:polyomino} and Bureaux and Enriquez \cite{Bureaux:polygons} refined this result by computing the asymptotic equivalent of the number of convex lattice chains. The former authors used analytic combinatorics tools similar to those we are using in this paper, and started from the study of the convex polyominoes. The latter started from Sinai's model, and used the local limit theorem from Bogatchev and Zarbaliev \cite{Bogatchev:limit}. They both got: 

$$p(n) \sim \frac{e^{- 2 \zeta'(-1)}}{(2\pi)^{7/6} \sqrt{3} \kappa^{1/18} n^{17/18}} \exp \left( 3 \kappa^{1/3} n^{2/3} + I\left(\left(\frac{\kappa}{n}\right)^{1/3}\right) \right), $$
with $I$ is a sum over the non-trivial zeros of the Riemann function $\zeta$ and under Riemann hypothesis, $I\left(\left(\frac{\kappa}{n}\right)^{1/3}\right)$ is of order $O(n^{1/6})$.\medskip


Eventually, computing the number of lattice polytopes still remains an open question as stepping up from convex lattice chains to convex polygons requires a thorough work over the convergence of the connections of the chains. \medskip

\subsection{Asymptotic estimate of the number of lattice zonotopes in $[0, n]^d$}
This paper is dealing with the enumeration of lattice zonotopes in $\R^d$, for any dimension $d$.  
A \emph{zonotope} is a convex geometric object defined as the Minkowski sum of $k$ segments, called its \emph{generators}. It is in one-to-one correspondence with the convex polygonal line made up by those generators. More particularly, a \emph{lattice zonotope} $Z$ (or integral zonotope) is a polytope for whom there exists $k \in \N$ and $(\vv_1, ..., \vv_k) \in (\Z^d)^k$ such that, up to translation, $Z$ is:

$$ Z = \left\{ \sum_{i=1}^k \alpha_i \vv_i  \:\:| \:\: \alpha_1,...,\alpha_k \in [0,1]\right\}. $$

As Minkowski sums of segments, zonotopes are more of a combinatorial object than a geometric one, and the question of enumerating them can be addressed. These objects are widely used in Combinatorics and other fields. Indeed they can be used when solving systems of polynomial equations \cite{Guibas:zonotopes} and for computing simplifications in spatial representation \cite{Huber:zonotopes}. From a theoretical perspective, they appear in the theory of the simplex method \cite{Delpia:simplex} and in the study of polytope diameters \cite{DezaPournin:diameter}. Evaluating their number is one more step in our comprehension of those objects and their relation to convex polytopes. \medskip

In 2018, B\'ar\'any, Bureaux, and Lund \cite[Theorem 1.2]{Barany:zonotopes} extended to higher dimensions the logarithmic asymptotic estimator of lattice zonotopes in a cone, ending at a given point. The authors computed the logarithmic equivalent of the number $z'_d(C, n \textbf{k})$ of lattice zonotopes inside a cone $C$ that starts at the origin of the cone and ends at a prescribed point $n\textbf{k}$. This gives the following result in the case of the zonotopes inscribed in a cube and ending at diagonally opposed points:  

$$ z'_d(\R_+^d, n \textbf{1}) \underset{n \rightarrow + \infty}{ = } \exp \left( (d+1) \kappa_d^{\frac{1}{d+1}}  n^{\frac{d}{d+1}} (1 + o(1)) \right),  \:\:\:\:\:\: \text{where} \:\: \kappa_d = \frac{\zeta(d+1)}{\zeta(d)}.$$

 Through a thorough analysis of the regularity of the generating function of $d$-dimensional lattice zonotopes, the present paper is devoted to refine the result of \cite{Barany:zonotopes}, and the result of \cite{Bureaux:polygons}, highlighting the polynomial that appears in the exponential starting at dimension 4. More precisely, we get: \medskip

\begin{thm}\label{theorem_1}
Denote the dimension $ d\in \N$, $d \geq 2$, and let $z_d(n)$ be the number of lattice zonotopes inscribed in $[0,n]^d$. 
There exist real numbers $\alpha_d$, $\beta_d$, $\kappa_d$ a polynomial $Q_d$ of degree $d$, and a function $I_{\text{crit}, d}$ depending on the non-trival zeros of the Riemann zeta function (all of which are explicitly given in the statement of Theorem \ref{theorem_1_complet}), such that, as $n$ grows large:

\begin{align}
     z_d(n) \sim  \alpha_d n ^{\beta_d} \exp \left( Q_d(n^\frac{1}{d+1}) + I_{\text{crit}, d}\left( \left(\frac{\kappa_d}{ n}\right)^{\frac{1}{d+1}} \right) \right),
\end{align}

\end{thm}

We hereby establish an asymptotic estimate in any given dimension.
Our work essentially deals with generating functions, broadly used in analytic combinatorics (see for instance Flajolet and Sedgewick \cite{Flajolet:AC}). This result naturally suggests two major remarks. 

\begin{rem}
Unexpectedly, the generalization of the exponential part is not a monomial as is in 2 and 3 dimensions (see below) but a polynomial of the $d+1$-th root of the size of the box. 
\end{rem}

\begin{rem}
As was already observed in dimension 2, the exact estimate contains a term depending on the set of non-trivial zeros of the $\zeta$ function. In Section \ref{section_numerical_discu}, we discuss the value of this term, and extend the results of \cite{Bodini:polyomino, Bureaux:polygons} that show that $I_{\text{crit}, 2}$ is negligible in dimension 2 for any "computably" large $n$ (typically  $< 10^{20}$).
\end{rem}

To illustrate this theorem, we give the following table:

\begin{table}[H]
    \centering
    \begin{tabular}{|c|c|c|c|}
        \hline
         & $\alpha_d$ & $\beta_d$ & $Q_d(X)$\\ [0.5ex] 
        \hline
        2D & $\frac{2^{1/9}  3^{13/18}  \zeta(3)^{2/9} e^{- 4 \zeta'(-1)} }{6 \pi^{16/9}}$  &  $- \frac{11}{9}$ & $\frac{2^{2/3}   3^{4/3}  \zeta(3)^{1/3} }{ \pi^{2/3}  }X^2 $  \\ [0.5ex] 
        \hline
        3D &  $\frac{2^{5/8}  3^{3/4}  5^{8/9} e^{- 4 \zeta'(-2)} }{120 \zeta(3)^{1/8} \pi}$ & $- \frac{13}{8}$ & $ \frac{2^{9/4}  \pi}{3^{1/2} 5^{1/3} \zeta(3)^{1/4}  }X^3$\\ [0.7ex] 
        \hline
        4D & $\left(\frac{2^{189} 3^{142} \zeta(5)^{71} }{5^{83}  \pi^{379}  } \right)^{\frac{1}{225}}  e^{- \frac{8}{3}\zeta'(-3) -\frac{16}{3} \zeta'(-1)}$ & $- \frac{521}{225}$ & $\frac{(5720)^{1/5} \zeta(5)^{1/5} }{\pi^{4/5}  } X^4 +  \frac{(16 720)^{2/5}  \zeta(5)^{2/5} \zeta(3) }{\pi^{18/5}  }X^2$\\
        \hline
    \end{tabular}
    \caption{Parameters of the asymptotic equivalent of the number of lattice zonotope in a hypercube of dimension 2, 3, and 4}
    \label{tab:number_of_zonotopes}
\end{table}
 
 
 

$I_{\text{crit}, 2}$, $I_{\text{crit},3}$, and $I_{\text{crit}, 4}$ are numerically approximated at $\left(\frac{\kappa_d}{ n}\right)^{\frac{1}{d+1}}$ in (\ref{numeric_I_crit}).  \medskip



\subsection{Moments of the diameter of lattice zonotopes.}
Once the combinatorial description of zonotopes and Theorem \ref{theorem_1} have been established, one can compute moments of combinatorial parameters of zonotopes. In this paper, we are interested in two parameters: the diameter of the graph of a lattice zonotope and the number of occurrences of a generator in a random lattice zonotope. \medskip

The diameter (in the sense of the diameter of the graph) of a polytope of a given size is a key combinatorial parameter, on which few things are known. Along with its own scientific interest (the famously now disproved Hirsch conjecture), it is connected to the complexity of the simplex algorithm. This makes upper bounds of this quantity actively looked for (see \cite{Delpia:simplex} for more references). In this perspective, zonotopes have been conjectured to be a class that reaches the largest possible diameter among all the lattice polytopes contained in $[0, n]^d$ (\cite[Conjecture 3.3]{Deza:primitive}). The largest possible diameter of a lattice zonotope contained $[0,n]^d$ is given in \cite{DezaPournin:primitivepoint}, while its exact asymptotic behavior is estimated in \cite{DezaPournin:diameter} when $d$ is fixed and $n \rightarrow + \infty$. \medskip

We hereby compute the asymptotic estimates of the mean of the distribution of the diameter of lattice zonotopes inscribed in $[0,n]^d$.

\begin{thm}\label{diameter_zono}
Let  $\mu_{diam}^n$ the mean of the distribution of the diameter of a lattice zonotope inscribed in $[0,n]^d$. Then, as $n$ grows large, we have:

$$ \mu_{diam}^n \underset{n \rightarrow +\infty}{= } \frac{\sqrt[d+1]{\kappa_d}}{\zeta(d+1)} \:\: n^{\frac{d}{d+1}} (1 + o(1)). $$
\end{thm}

This article is organized as follows. We begin, in Section \ref{section_combi}, with establishing the generating function associated to the combinatorial class of lattice zonotopes. Theorem \ref{theorem_1} is established in Section \ref{section_proof_main} using multidimensional saddle-point method. In order to use this method, Section \ref{section_analysis} is dedicated to find an equivalent to the generating function in the univariate case, and to prove some technical lemmas to control the derivative of the generating function.
Finally, the last two sections are respectively dedicated to discussing the obtained asymptotic equivalent and to computing the estimated moments of the number of non-colinear generators and the number of lattice points crossed by a given generator.   




\section{CLASS OF LATTICE ZONOTOPES, COMBINATORIAL APPROACH}\label{section_combi}
Remind from the introduction that a \emph{lattice zonotope} $Z$ (or integral zonotope) is a polytope for whom there exists $k \in \N$ and $\vv_1, ..., \vv_k \in \Z^d$ such that, up to translation, $Z$ is:

$$ Z = \left\{ \sum_{i=1}^k \alpha_i \vv_i  \:\:| \:\: \alpha_1,...,\alpha_k \in [0,1]\right\}. $$

As we study inscribed zonotopes in a hypercube, we consider lattice zonotopes up to translation thereafter.  That is to say that all the lattice segments start at $0$.

\begin{figure}[H]
    \centering
    \subfloat[\centering 2 dimensions]{{\includegraphics[width=5cm]{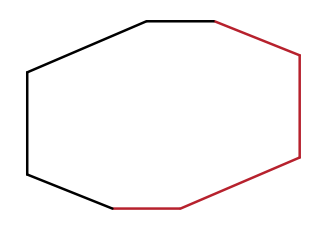} }}%
    \qquad
    \subfloat[\centering 3 dimensions]{{\includegraphics[width=5cm]{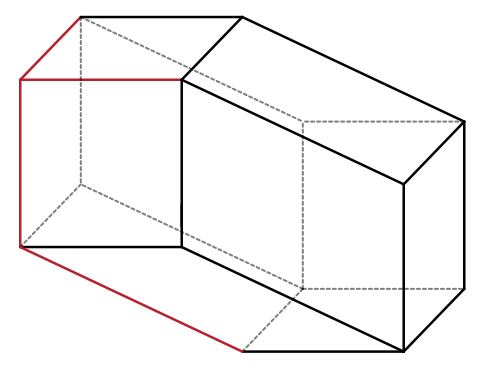} }}%
    \caption{Zonotopes made from 4 red-colored generators.}%
    \label{fig:example}%
\end{figure}

Given a lattice zonotope $Z$ and generators $\vv_1,..,\vv_k \in \Z^d\setminus\{\boldsymbol{0}\}$, the set of generators uniquely defines a zonotope but the inverse is not true. We denote $\textsc{Zon}_d$ the set of $d$ dimensional lattice zonotopes up to translation and will use the $L^1$ norm of its generators, $\sum_{i=0}^k |\vv_i|$, as the size of a lattice zonotope $Z$.
As it was previously established in study of polygonal chains \cite{Bureaux:polygons}, we state a bijective correspondence between zonotopes and finite support functions encoding the set of generators. \medskip

A vector $\vv$ is named a primitive vector if the gcd of its non-zero coordinates is 1, namely if the intersection between the segment $[0,\vv]$ and $\Z^d$ is exactly the extremities 0 and $\vv$. Let $\P_d$, respectively $\P_{d\: +}$ denote the set of primitive vectors of $\Z^d$, respectively $(\Z_+)^d$. Then let $\mathcal{P}_d\star$ denote the set of primitive vectors whose first non-zero coordinate is positive, following the convention from \cite{Deza:primitive}. The space $\Omega$ of non-negative integer-valued functions $\omega : \mathcal{P}_d\star \rightarrow \N$ with finite support is in one-to-one correspondence with $d$ dimensional lattice zonotopes up to translation $\textsc{Zon}_d$. To be more specific:\medskip

\begin{itemize}
    \item The function $\omega$ associated to a zonotope $Z$ is defined for all $\vv \in \mathcal{P}_d\star$ as the number of occurrence of $\vv$ in an edge of $Z$ that is colinear to $\vv$ if there is such an edge, and 0 otherwise. 
    
    \item  Given a function $\omega$, the list of generator which defines $Z$ is $\{ \omega(\vv) \vv, \vv \in \text{supp}(\omega)\}$. In particular, $\underset{\vv\in  \mathcal{P}_d\star}{\sum}\omega(\vv)\vv$ is the vertex of the zonotope opposite the origin (which we will refer to as the endpoint of the zonotope).
\end{itemize}

This combinatorial structure can be summed up in a set of repeated elements without taking care of the order. This coincides with the notion of a Multiset (defined in \cite[Chapter~1]{Flajolet:AC}, p.26-27) in the framework of the symbolic method. This gives us the possibility, in this paper, to apply well known analytical tools on the generating function of the multiset of elements of  $ \mathcal{P}_d\star$, the study of which is a part of the contribution of this paper. \medskip

Let $Zon_d$ encode the multivariate generating function of the class $\textsc{Zon}_d$. We will show in (\ref{partition_function}) that this function is defined on the open centered disk of radius $1$, hence for $\xx = (x_1,x_2,..., x_d) \in (-1,1)^d$, the generating function is defined as

\begin{align}\label{formula_defin_z_n}
     Zon_d(\xx) = \sum_{\nn \in \N^d} z_{\nn} x_1 ^{n_1}... x_d^{n_d},
\end{align}

where $z_{\nn}$ is the number of lattice zonotopes inscribed in a box of size $n_1 \times n_2 \times ... \times n_d$. The variable $x_i$ then encodes the length of the projection of a zonotope on the $i$th axis of coordinates. We classically denote $\textbf{\textit{x}}^\textbf{\textbf{n}} =x_1 ^{n_1}... x_d^{n_d}$, $\un$ the vector $(1,1, ...,1)$, and $\xx \cdot \yy$ as the canonical scalar product of $\xx$ and $\yy$ in the sequel. We also introduce the notation $[\vv]$ to designate the vector of the absolute values of $\vv$, $(|v_1|,..., |v_d|)$. \medskip

Adding a given generator $\omega(\vv)\vv$ to a zonotope increases the size of the box containing the zonotope by $\omega(\vv)|v_1| \times ... \times \omega(\vv)|v_d|$. Then we can rewrite this functions using the finite support function as follows:

$$Zon_d(\xx) =  \sum_{\omega \in \Omega} \prod_{\vv \in \mathcal{P}_d\star} \xx^{\omega(\vv) [\vv]}, $$

 which can be transformed into the generating function of a multiset by factorization and elementary operations:
 
\begin{align}\label{partition_function}
 Zon_d(\xx) =  \prod_{ \vv \in \mathcal{P}_d\star}\left(1 - \xx^{[\vv]} \right)^{-1}.   
\end{align}

This function is defined on the open disk of radius 1 and centered at 0, with a singularity at $1$. In the following, we use the change of variables $x_i = e^{-\theta_i}$ to study the function as $\ttheta$ goes to 0. The resulting form is a partition function over the set $\mathcal{P}_d\star$ of primitive vectors whose first non-zero coordinate is positive. This function is analogous to the partition function used in the Boltzmann probabilistic point of view in \cite{Bureaux:polygons}, that is, for $\ttheta \in (0, +\infty)^d$:

$$ \prod_{ \vv \in \mathcal{P}_d\star}\left(1 - e^{ - \ttheta \cdot [\vv]} \right)^{-1}. $$

Denoting $d(\vv)$ the number of non-zero coordinates of $\vv$, the product can be simplified by gathering in the same factor the vectors whose absolute value of coordinates coincide (the first non-zero coordinate shall be positive), that is 

\begin{align}\label{generating_function}
    Zon_d\left(e^{- \ttheta}\right) =  \prod_{ \vv\in \P_{d+}} \left( 1 - e^{-  \ttheta \cdot \vv } \right)^{- 2^{d(\vv)-1}}.
\end{align}




\section{ASYMPTOTIC ANALYSIS OF THE GENERATING FUNCTION}\label{section_analysis}
We start the analysis of the generating function by giving an integral formula, and an equivalent of the logarithm of $Zon_d$ and of its partial derivative. The second part focuses on the study of the univariate generating function defined by $Zon_d(x\boldsymbol{1})$. We compute the asymptotic equivalent of this function. Those 2 parts are fundamental in the convergence analysis conducted in the next section. 

\subsection{Integral formula and equivalent of partial derivative of the generating function.}
In order to asymptotically study the generating function, we will mainly use in the following an integral formula, deriving from the Mellin inversion formula:

\begin{lem}\label{lemma_integral_form_gg}
Taking $c>d$, for all $\ttheta \in (0, + \infty)^d$, we have
\begin{equation}\label{equa_integral}
  \log\left(  Zon_d(e^{-\ttheta}) \right) =   \frac{1}{2 i \pi} \int _{c - i\infty} ^{ c+i \infty} \frac{\zeta(s+1)\Gamma(s)}{\zeta(s)}  \sum_{ \vv\in \Z^d_+ \setminus \{\boldsymbol{0}\}} \frac{2^{d(\vv)-1}}{( \ttheta \cdot \vv)^s}  ds.  
\end{equation}

\end{lem}

\begin{proof}
Let $\ttheta \in (0, + \infty)^d$, we first compute the logarithm of the equation (\ref{generating_function}), depending on the number $d(\xx)$ of non-zero coordinate of $\vv$, and write the Taylor series expansion of the logarithm:

$$ \log\left(  Zon_d(e^{-\ttheta}) \right) = - \sum_{ \vv\in \P_{d+}} 2^{d(\vv)-1} \log \left( 1 - e^{-  \ttheta \cdot \vv } \right)  = \sum_{ \vv\in \P_{d+}} \sum_{m \geq 1} \frac{ 2^{d(\vv)-1} }{m} e^{-m  \ttheta \cdot \vv }   $$

Now, we can compute the Mellin transform of this function by replacing the exponential terms with their transform through Mellin inversion formula. Recall that the Mellin transform of the real function $y \mapsto e^{-y}$ is the $\Gamma$ function, hence by the inverse of the Mellin transform, for every real numbers $y > 0$ and $c>0$:

$$e^{-y} =  \frac{1}{2 i \pi} \int _{c - i\infty} ^{ c+i \infty} \Gamma(s) y^{-s} ds.$$

We replace each $e^{-m  \ttheta \cdot \vv }$ by using this relation. Then, for any real number $c$ such that $c > d$, we have $\sup_{s \in c + i \R}\left(\sum_{\vv \in  \P_{d+} } \frac{1}{(\ttheta \cdot \vv )^{s}}\right) <+ \infty$ and $ \sup_{s \in c + i\R} \left( \sum_{ m\geq1} \frac{1}{m^{s+1}}\right) < + \infty$. Therefore the Fubini-Tonelli theorem leads to:

$$ \log\left(  Zon_d(e^{-\ttheta}) \right) =   \frac{1}{2 i \pi} \int _{c - i\infty} ^{ c+i \infty}  \sum_{m \geq 1}  \sum_{ \vv\in \P_{d+}} \frac{ 2^{d(\vv)-1} }{m} \frac{\Gamma(s)}{(m \ttheta \cdot \vv)^s}  ds. $$

We can extract the variable $m$ from the sums and bring out the Riemann $\zeta$ function. We obtain 

$$ \log\left(  Zon_d(e^{-\ttheta}) \right) =   \frac{1}{2 i \pi} \int _{c - i\infty} ^{ c+i \infty}  \zeta(s+1)  \sum_{ \vv\in \P_{d+}} \frac{ 2^{d(\vv)-1} \Gamma(s)}{( \ttheta \cdot \vv)^s}  ds. $$

Finally, the sum over the primitive vectors can be completed to a sum over $\Z^d_+$ the using the partition $\Z^d_+ = \cup_{k \geq 1} k \P_{d+}$,

\begin{align}\label{partitionZd}
    \sum_{ \vv\in \Z^d_+ \setminus \{\boldsymbol{0}\}} \frac{2^{d(\vv)-1}}{( \ttheta \cdot \vv)^s} = \sum_{k \geq 1 } \frac{1}{k^s} \sum_{ \vv\in \P_{d+} } \frac{2^{d(\vv)-1}}{( \ttheta \cdot \vv)^s} = \zeta(s) \sum_{ \vv\in \P_{d+} } \frac{2^{d(\vv)-1}}{( \ttheta \cdot \vv)^s}. 
\end{align}

For $\ttheta \in (0, + \infty)^d$, the integral becomes, after introducing a term $\zeta(s)$ to use (\ref{partitionZd}), 

$$ \log\left(  Zon_d(e^{-\ttheta}) \right) =   \frac{1}{2 i \pi} \int _{c - i\infty} ^{ c+i \infty} \frac{\zeta(s+1)\Gamma(s)}{\zeta(s)}  \sum_{ \vv\in \Z^d_+ \setminus \{\boldsymbol{0}\}} \frac{2^{d(\vv)-1}}{( \ttheta \cdot \vv)^s}  ds.  $$

\end{proof}

The following lemma establishes the asymptotic equivalent of the partial derivatives of $\log\left( Zon_d(e^{-\ttheta}) \right)$. These asymptotic equivalents of partial derivative are crucial in the proof of Proposition \ref{definition_gittenberger} that gives an asymptotic equivalent of the coefficients of $Zon_d$.

\begin{lem}\label{lemma_derivativ}
For $(k_1,...,k_d) \in \N^d$ and for all $\epsilon > 1$, for $\ttheta \in (0, +\infty)^d$ such that for $1\leq i, j \leq d$ and $i\neq j$, $\frac{\theta_i}{\theta_j} \in (\frac{1}{\epsilon}, \epsilon)$,  

$$\frac{\partial^{k_1 + k_2 +...+ k_d}}{\partial\theta_1^{k_1} \partial\theta_2^{k_2}...\partial \theta_d^{k_d}} \log\left(  Zon_d(e^{-\ttheta}) \right) \underset{\ttheta \rightarrow \boldsymbol{0}}{\sim} (-1)^{k_1 + ... + k_d} \frac{2^{d-1}\zeta(d+1)}{\zeta(d)}\frac{\prod_{k=1}^d k_i !}{\theta_1^{k_1+1}...\theta_d^{k_d+1}}.$$
\end{lem}

\begin{proof}
We introduce the Barnes zeta function in $d$ dimensions $\zeta_d(s, w, \ttheta)$, where $s$, $w$ and  $\{\theta_i, 1\leq i\leq d \}$ are complex numbers such that $\Re(s) > d$, $\Re(w) >0$ and for $1 \leq i \leq d$, $\Re(\theta_i)> 0$ as:

$$\zeta_d(s, w, \ttheta) = \sum_{\vv \in \Z_+^d} \frac{1}{(w + \ttheta \cdot \vv)^s}. $$

This function can be meromorphically continued to all complex $s$, with simple poles at $1, 2, ..., d$. Remark that all the vectors $\vv$ for which $d(\vv) = d$ are the vectors without any zero coordinate, that is vectors belonging to $\Z_+^d \cap (0,+\infty)^d$. Additionally we can write

$$
    \sum_{\vv \in \Z_+^d \cap (0,+\infty)^d} \frac{2^{d(\vv)-1}}{(\ttheta \cdot \vv)^s} = 2^{d-1} \zeta_d\left(s, \sum_{i= 1}^d \theta_i, \ttheta \right ).
$$

This leads to rewriting the sum in the integrand of the right term of (\ref{equa_integral}) in Lemma \ref{lemma_integral_form_gg} as a sum of Barnes zeta functions. For $s > d$, that is:

$$ \sum_{ \vv\in \Z^d_+ \setminus \{\boldsymbol{0}\}} \frac{2^{d(\vv)-1}}{( \ttheta \cdot \vv)^s} = \sum_{j = 1}^d \sum_{1 \leq i_1 <... <i_j \leq d} 2^{j -1} \zeta_i \left(s, \sum_{k=1}^j \theta_{i_k} , (\theta_{i_1},...,\theta_{i_2}) \right).$$

The Barnes $\zeta$ function in $\delta$ dimensions has a simple pole at $s \in \{1, 2, ...,  \delta \}$, hence the only pole in the integrand comes from the $d$-dimensional Barnes $\zeta$ function. We translate the domain of integration $(c-i\infty,c+i\infty)$ by 1 to the left, by using the residue theorem at the pole at $s = d$. For $d-1 < c < d$, we have: 

$$ \log\left(  Zon_d(e^{-\ttheta}) \right) = \frac{ 2^{d-1} \zeta(d+1)\Gamma(d)}{\zeta(d) \prod_{i=1}^d \theta_i}  + \frac{1}{2 i \pi} \int _{c - i\infty} ^{ c+i \infty} \frac{\zeta(s+1)\Gamma(s)}{\zeta(s)}  \sum_{ \vv\in \Z^d_+ \setminus \{\boldsymbol{0}\}} \frac{2^{d(\vv)-1}}{( \ttheta \cdot \vv)^s} ds.  $$

With the bound established in Lemma \ref{lemma_Barnes} (below), the integral in the right term of this equation is differentiable, hence the right term is differentiable, and therefore we can differentiate it. The derivative of the integral is negligible in front of the leading term, so we get the wanted result. 

\end{proof}

\begin{lem}\label{lemma_Barnes}
For all $\kk \in \N^d$, $\delta \in (0,1)$, and $\epsilon > 0$, there exists a constant $C> 0$ such that all $\ttheta$ with their coordinates $(\theta_i)$ respecting $\epsilon < \theta_i < \frac{1}{\epsilon}$ satisfy, for all $s$ such that $\Re(s) \in (d-1 + \delta, d+1)$

$$ \left| \frac{\partial^{\kk \cdot \boldsymbol{1}}}{ (\boldsymbol{\partial} \ttheta)^{\kk} } \zeta_{d}\left( s, \ttheta\cdot \boldsymbol{1}, \ttheta \right) \right| \leq \frac{C |s|^C}{|\ttheta|^{\kk \cdot \boldsymbol{1} + \Re(s)} |s-d|}.$$
\end{lem}

\begin{proof}
Let $\{x\}= x - \lfloor x \rfloor $ denote the fractional part of $x$. For $x\in \R_{+}$, consider the function defined by $ x \mapsto F(x)= \sum_{n_2, ...,n_d \geq 0} (w + \theta_1 x + \theta_2 n_2 + ... + \theta_d n_d)^{-s}$. We apply the Euler--Maclaurin formula to the $d$-dimensional Barnes $\zeta$ function using $F$, leading to

$$\zeta_d(s,w, \ttheta) = \sum_{n_1 \geq 0} F(n_1) = \int_0^{\infty} F(x)dx - \frac{F(0)}{2} + \int_0^{\infty}\left(\{x\} - \frac{1}{2}\right)F'(x) dx. $$

We rewrite the right term as the following expression, because $F(x) = \zeta_{d-1}(s, w+ \theta_1 x, (\theta_2, ...\theta_d))$, for $\Re(w) >0$ and $\Re(\ttheta_i)> 0$. For each $1 \leq i \leq d$:

\begin{align}\label{Barnes_d_dim}
    \zeta_d(s,w, \ttheta)  = \int_0^{\infty} \zeta_{d-1}(s, &w+ \theta_1 x, (\theta_2, ...\theta_d))dx - \frac{\zeta_{d-1}(s,w, (\theta_2, ...\theta_d))}{2} \\
    & -s \theta_1 \int_0^{\infty} \left( \{ x \} - \frac{1}{2}\right) \zeta_{d-1}(s+1, w+ \theta_1 x, (\theta_2, ...\theta_d)) dx. \nonumber
\end{align}

As mentioned before, the $d$-dimensional Barnes $\zeta$ function of parameters $s$, $w$, and $\ttheta$ has a meromorphic continuation to all complex $s$ whose only singularities are simple poles at $1, 2, ..., \text{ and }d$. This implies that only the first term in the previous expression has a pole at $s = d$ while the other two terms have poles at $s =d-1, d-2, ..., 1$. \medskip

We can use again Euler--Maclaurin formula for each of $\zeta_{d-1}$. We recall the case $d=2$, which is given in the Lemma A.1 in \cite{Bureaux:polygons}:

\begin{align}\label{Barnes_deux_dim}
    \zeta_2(s, \omega, (\theta_1, \theta_2)) =& \frac{1}{\theta_1 \theta_2}\frac{\omega^{-s+2}}{(s-1)(s-2)} + \frac{(\theta_1 + \theta_2)\omega^{-s+1}}{2\theta_1 \theta_2(s-1)} + \frac{\omega^{-s}}{4} \nonumber\\
    &- \frac{\theta_2}{\theta_1}\int_{0}^{+\infty}\frac{\{y\} -\frac{1}{2}}{\omega+ \theta_2 y)^s}dy - \frac{\theta_1}{\theta_2} \int_{0}^{+\infty}\frac{\{x\} -\frac{1}{2}}{\omega+ \theta_1 x)^s}dx \\
    & - s \frac{\theta_2}{2} \int_{0}^{+\infty}\frac{\{y\} -\frac{1}{2}}{\omega+ \theta_2y)^{s+1}}dy - s \frac{\theta_1}{2} \int_{0}^{+\infty}\frac{\{x\} -\frac{1}{2}}{(\omega+ \theta_1 x)^{s+1}}dx \nonumber\\
    & + s(s+1)\theta_1\theta_2 \int_{0}^{+\infty} \int_{0}^{+\infty}\frac{(\{x\} -\frac{1}{2})(\{y\} -\frac{1}{2})}{(\omega + \theta_1 x + \theta_2 y)^{s+2}}dx dy\nonumber
\end{align}

Therefore, with recursive use of the Euler--Maclaurin formula on the each of the three terms of (\ref{Barnes_d_dim}) and so after, we obtain a formula with a finite number of terms, linearly depending on the $(\theta_i)$. The first term is

$$ \frac{1}{\theta_1...\theta_d} \frac{w^{-s+d}}{(s-1)...(s-d)}. $$

All the other terms of the development are of the form of the terms of (\ref{Barnes_deux_dim}), with quotients of $\theta_i$ and $\theta_j$. All these terms are differentiable, and therefore so is the right term of (\ref{Barnes_deux_dim}). By assumptions, all ratios  $\theta_i/ \theta_j$ range between $\epsilon$ and $\frac{1}{\epsilon}$, hence we can find a constant $C$ such that all terms are upper bounded by 

$$\frac{C |s|^C}{|\ttheta|^{\kk \cdot \boldsymbol{1} + \Re(s)} |s-d|}. $$

\end{proof}

\subsection{Asymptotic equivalent of the univariate generating function.}\label{subsection_univariable_equival}

We now can compute the asymptotic equivalent of the generating function when all variables $\theta_i$ are equal and tend to 1. It will be used in the proof of the main theorem to compute the equivalent of the coefficient of the generating function given in Proposition \ref{definition_gittenberger}. In the process, we introduce the following notation:

\begin{defn}\label{definition_operator}
Let $\mathrm{M}(\C)$ be the field of meromorphic functions in $\C$. We define the polynomial $P_d({X})$ as

$$P_d({X}) = \sum_{\delta = 1}^d  \binom{d}{\delta} 2^{\delta-1} \frac{\prod_{k =1}^{\delta-1}(X -k)}{(\delta-1) !} = p_{d,d-1} {X}^{d-1} + ... + p_{d,1} {X} + p_{d,0}, $$ 

with the convention that $\prod_{k =1}^{\delta-1}(X -k) = 1 $ if $\delta =1$, and the operator $\mathbf{\Pi}_d$ as
\begin{align*}
   \mathbf{\Pi}_d \colon \mathrm{M}(\C) &\to \mathrm{M}(\C)\\
  \phi &\mapsto p_{d,d-1} \phi(\cdot \: - (d-1)) + ... + p_{d,1} \phi(\cdot \: - 1) + p_{d,0} \phi.
\end{align*}
\end{defn}

\begin{rem}
The family of polynomials $(P_d(X) )_{d\geq 1}$ is recursively defined by $P_{d+2} (X) = \frac{2X}{d+1}P_{d+1}(X) + P_d(X)$, with $P_1({X}) = 1$  and $P_2({X}) = 2 {X}$. Which results in alternating odd and even polynomials, the even ones having 1 as constant term. Moreover, the leading term of $P_d$ is of degree $d-1$ and its coefficient is $\frac{2^{d-1}}{(d-1)!}$.
\end{rem}

\begin{table}[H]
    \centering
    \begin{tabular}{|c|c| c|c|}
        \hline
        $P_1$ & $1$ & $\Pi_1(\phi )$ & $\phi$ \\ 
        \hline
        $P_2$ & $2X$ & $\Pi_2(\phi )$ & $2 \phi( \cdot - 1)$\\ 
        \hline
        $P_3$ &  $\frac{4}{3} X^2 + 1 $ & $\Pi_3(\phi )$ & $\frac{4}{3} \phi( \cdot - 2) + \phi $\\ 
        \hline
        $P_4$ & $ \frac{4}{3} X^3 + \frac{8}{3} X$ & $\Pi_4(\phi )$ &  $ \frac{4}{3} \phi( \cdot - 3) + \frac{8}{3} \phi( \cdot - 1)$\\
        \hline
    \end{tabular}
    \caption{First polynomials $P_d$.}
    \label{tab:polynomials}
\end{table}

\begin{figure}[H]
    \centering
    \begin{overpic}[width = 13cm]{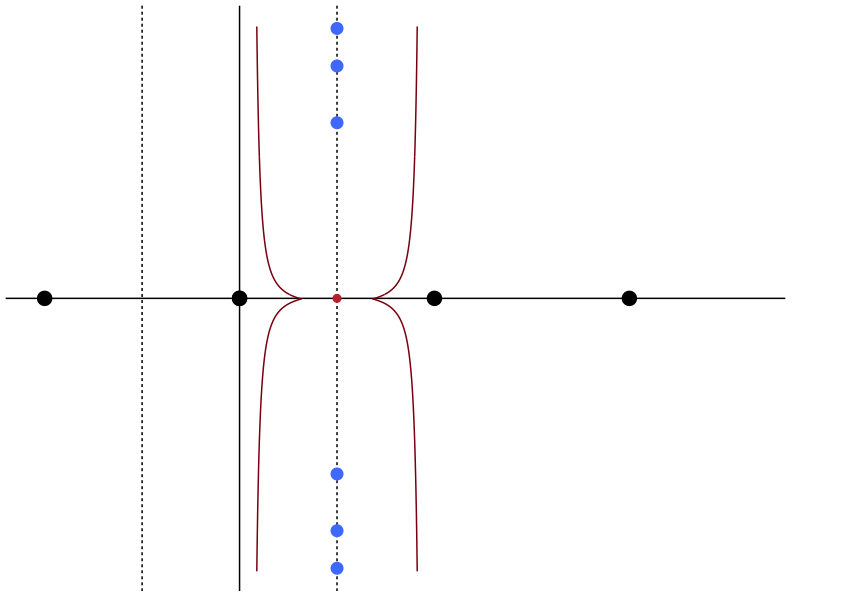}
        \put(25,31){0}
        \put(49,31){1}
        \put(72,31){2}
        \put(-1,31){$-1$}
        \put(37,31){$\frac{1}{2}$}
        \put(10,31){$-\frac{1}{2}$}
        \put(50,10){\color{BrickRed}$\gamma_{r}$}
        \put(31,10){\color{BrickRed}$\gamma_{l}$}
        \put(70.7,6){\color{Apricot}\rule{1pt}{217.1pt}}
        \put(16.8,6){\color{Apricot}\rule{200pt}{1pt}}
        \put(16.8,6){\color{Apricot}\rule{1pt}{217.1pt}}
        \put(16.8,64.5){\color{Apricot}\rule{200pt}{1pt}}
        \put(0,6){\color{Apricot}$-\frac{1}{2}- iT_{k}$}
        \put(0,64.5){\color{Apricot}$-\frac{1}{2}+ iT_{k}$}
    \end{overpic}
    \caption{The contour $\gamma$ surrounding the critical strip. The blue points stand for the non trivial zeros of the Riemann $\zeta$ function.}
    \label{fig:bande_critique}
\end{figure}
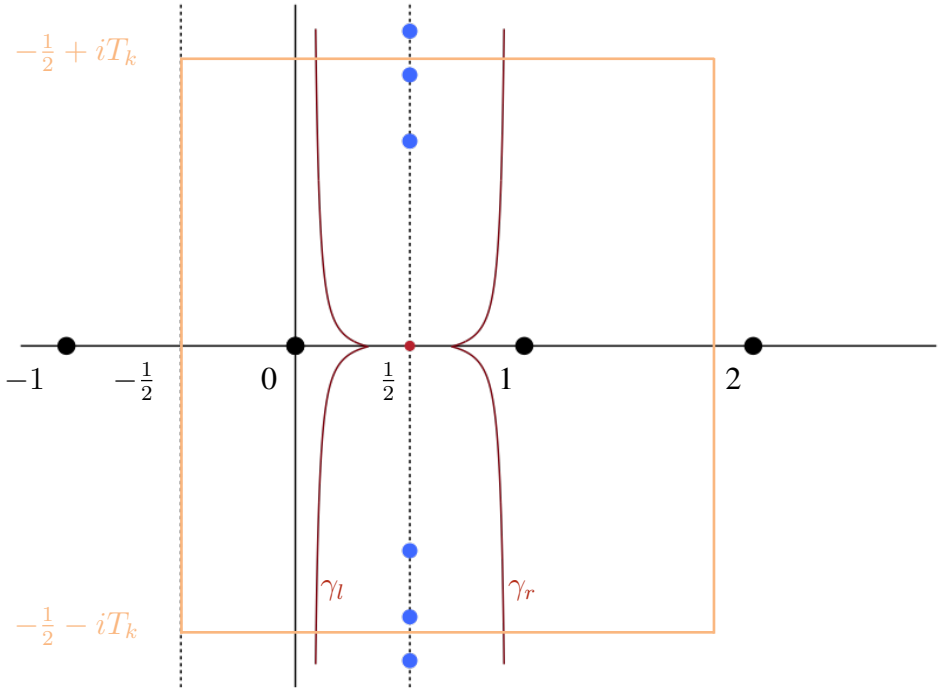

The operator naturally results from the Mellin transform, as for $s\neq 1$, $\mathbf{\Pi}_d (\zeta)(s) = \sum_{k \geq 1} \frac{P_d(k)}{k^s}$.
With this notation, we present the following proposition, where the polynomial in $\frac{1}{\theta}$ is the origin of the polynomial in the exponential part of Theorem \ref{theorem_1}:

\begin{prop} \label{prop_equivalent_univariate}
Let $\gamma$ denote the contour defined as the union of the two oriented paths, depending on an integer $A>0$ chosen such that it surrounds all the zeros of the Riemann $\zeta$ function : on the right  the curve $\gamma_{r}(t) = 1 - \frac{A}{\log(2 + |t|)} + it $ for $t$ going from $- \infty$ to $+\infty$, and on the left the curve $\gamma_{l}(t) = \frac{A}{\log(2 + |t|)} + it $ for $t$ going from $+ \infty$ to $- \infty$ (see Figure \ref{fig:bande_critique}).\medskip

We define the function $ I_{\text{crit}, d}$ by 

\begin{align}
    I_{\text{crit}, d}(\theta) = \frac{1}{2 i \pi } \int_{\gamma}  \frac{\mathbf{\Pi}_d (\zeta(s)) }{ \zeta(s)} \zeta(s+1) \Gamma(s) \theta^{-s} ds 
\end{align}
For all $\theta> 0$, the asymptotic equivalent of $ Zon\left(e^{-\theta}\right)$ when $\theta\rightarrow 0$ is:

\begin{align}\label{expr:prop:estimation_logZ}
    \log \left(Zon \left(e^{-\theta} \right) \right) = \sum_{\delta=1}^{d-1}\left(p_{d,\delta}\frac{\zeta(\delta+2) \Gamma(\delta+1)}{\zeta(\delta+1) \theta^{\delta+1}} \right) + I_{\text{crit}, d}(\theta) +  C + O \left( n^{- \frac{1}{d+1}}\right) 
\end{align}

with $C = 2 \mathbf{\Pi}_d \left(\log(2\pi) \zeta -  \zeta'\right) (0) + 2 \mathbf{\Pi}_d \left(\zeta\right)(0) \log(\theta) $.\medskip

\end{prop}

\begin{proof}
For any real number $c$ such that $c> d$, and $\theta \in (0, +\infty)$, we fix $\ttheta = \theta \boldsymbol{1}$ and express the integral expression (\ref{equa_integral}) in Lemma \ref{lemma_integral_form_gg}  as:

$$  \log\left(  Zon_d(e^{-\theta} \boldsymbol{1}) \right) =   \frac{1}{2 i \pi} \int _{c - i\infty} ^{ c+i \infty} \frac{\zeta(s+1)\Gamma(s)}{\zeta(s)}  \sum_{ \vv\in \Z^d_+ \setminus \{\boldsymbol{0}\}} \frac{2^{d(\vv)-1}}{( \theta \boldsymbol{1} \cdot \vv)^s}  ds.  $$

The vector $\vv$ can be clustered with the $\binom{d}{d(\vv)}-1$ other vectors with the same non-zero ordered coordinates, as the scalar product $\boldsymbol{1}\cdot \vv$ is equal for all this vectors. We deduce that

$$  \sum_{ \vv\in \Z^d_+ \setminus \{\boldsymbol{0}\}} \frac{2^{d(\vv)-1}}{( \boldsymbol{1} \cdot \vv)^s} = \sum_{\delta = 1}^{d} 2^{\delta-1} \binom{d}{\delta}  \sum_{ \vv\in \Z^{\delta}_+ \cap (0, + \infty )^{\delta} } \frac{1}{( \boldsymbol{1} \cdot \vv)^s}.$$

Then, we can turn the sum over $\Z^{\delta}_+ \cap (0, + \infty )^{\delta} $ into a sum labelled by the scalar product value, denoted $n \in \N$:

$$   \sum_{ \vv\in \Z^{\delta}_+ \cap (0, + \infty )^{\delta} } \frac{1}{( \boldsymbol{1} \cdot \vv)^s} =  \sum_{ n \geq 1} \frac{\binom{\delta-1}{n-1}}{ n^s}.$$

Hence, with the operator $\mathbf{\Pi}_d$ defined in Definition \ref{definition_operator}, we can write, for $c > d$:

\begin{align}\label{integral_expression_univar}
      \log\left(  Zon_d(e^{-\theta} \boldsymbol{1}) \right) =   \frac{1}{2 i \pi} \int _{c - i\infty} ^{ c+i \infty} \frac{\mathbf{\Pi}_d (\zeta(s)) }{ \zeta(s) \theta^s} \zeta(s+1) \Gamma(s) ds, 
\end{align}

We observe that (\ref{integral_expression_univar}) is the expression of an inverse Mellin transform. In other words, the Mellin transform of the left hand side of (\ref{integral_expression_univar}) is $\theta^s$ times the integrand in the right-hand side: for $s\in \C$, and $\Re(s) > d$, 

$$\mathcal{M}\left[\log \left(Zon_d \left(e^{-\theta} \un \right)\right)\right] (s)  = \int_0^{+\infty} \log \left(Zon_d \left(e^{-\theta} \un \right)\right)\theta^{s-1}d\theta =\frac{\mathbf{\Pi}_d (\zeta(s)) }{ \zeta(s)} \zeta(s+1) \Gamma(s). $$

In order to give some concrete idea of this formula, here are the Mellin transform of the 2, 3, and 4 dimensional cases:

\begin{align*}
    \mathcal{M}\left[\log \left(Zon_2 \left(e^{-\theta} \un\right)\right)\right] (s) &=\frac{2 \zeta(s - 1) }{ \zeta(s)} \zeta(s+1) \Gamma(s),\\
    \mathcal{M}\left[\log \left(Zon_3 \left(e^{-\theta} \un \right)\right)\right] (s) &=\frac{2 \zeta(s - 2) + \zeta(s) }{ \zeta(s)} \zeta(s+1) \Gamma(s),\\
    \mathcal{M}\left[\log \left(Zon_4 \left(e^{-\theta} \un \right)\right)\right] (s) &=\frac{ \frac{4}{3}\zeta(s -3) + \frac{8}{3}\zeta(s -1) }{ \zeta(s)} \zeta(s+1) \Gamma(s).
\end{align*}

 All the function $\zeta$, $\Gamma$, and $1/\zeta$ that compose $ \mathcal{M}\left[\log \left(Zon_d \left(e^{-\theta} \un \right)\right)\right]$  can be continued into meromorphic ones on $\C$, and so do the Mellin transform. To obtain the successive orders of the right term of (\ref{expr:prop:estimation_logZ}), we shift the vertical line to the left (in the complex plane), and use the residue theorems on the poles of the integrand of (\ref{integral_expression_univar}) to switch the line of integration from the right to the left each pole of the Mellin transform. This method is widely known as utmost-left propagation of the integration contour in a transform inversion formula with the residue theorem around the poles (see \cite[p. 765]{Flajolet:AC}). The poles of $ \mathcal{M}\left[\log \left(Zon_d \left(e^{-\theta} \un \right)\right)\right]$ are located at each real number $\delta$, $\delta$ being an integer between 0 and $d$, and at Riemann's $\zeta$ function's non-trivial zeros due to the denominator. \medskip

The integral of $ \mathcal{M}\left[\log \left(Zon_d \left(e^{-\theta} \un \right)\right)\right] \theta^{-s}$ is well defined on the vertical lines with real part in $(i, i+1)$ for $i\in \llbracket 1 , d-1 \rrbracket$ because of exponential decrease of $\Gamma(s)$ when $\Im(s)$ goes to $\pm \infty$. Starting from (\ref{integral_expression_univar}), we use the residue theorem around the pole at $d$ to shift the line of integration of equation $\Re(s) = c$, with $c >d$, to the line of equation $\Re(s) = c_1$, with $c_1\in (d-1, d)$. It follows that, for this $c_1\in (d-1, d)$:

$$ \log(Zon_d (e^{-\theta} \un )) = \frac{p_{d, d-1} \zeta(d +1) \Gamma(d) }{\zeta(d) \theta^{d}} + \frac{1}{2i\pi} \int_{c_1- i\infty}^{c_1+ i\infty} \frac{\mathbf{\Pi}_d (\zeta(s)) }{ \zeta(s)} \zeta(s+1) \Gamma(s) \theta^{-s} d s.$$

We can repeat the process until $1 < c_2 <2$. The last step is analogous to the proof of Lemma 2.2 in \cite{Bureaux:polygons}. 
Recall that $\gamma$ is the contour defined in Proposition \ref{prop_equivalent_univariate}, parametered by $A$. The existence of such $A$ is proven in \cite[Theorem 3.8]{Titchmarsh:riemann} (we only need the existence of $A$ here, as we don't need it when we discuss the value of $I_{\text{crit}, d}$ in Section \ref{section_numerical_discu}). \medskip

Hereafter, we prove that the integral of $ \mathcal{M}\left[\log \left(Zon_d \left(e^{-\theta} \un \right)\right)\right]$ on the line of equation $\Re(s) = c_2$ is the sum of the residue of the integrand around $s = 1$, $I_{\text{crit}, d}$, and the integral of $ \mathcal{M}\left[\log \left(Zon_d \left(e^{-\theta} \un \right)\right)\right]$  on the line of equation $\Re(s) = -\frac{1}{2}$.\medskip 

The latter integral is convergent because $\Gamma(s)$ exponentially decreases when $\Im(s) \rightarrow \pm \infty$ and $1/ \zeta(s)$ is bounded on the line of equation $\Re(s) = -\frac{1}{2}$. With this bound, the integrand divided by $\theta^{\frac{1}{2}}$ is dominated, and we obtain that the integral on the domain of equation $\Re(s)= -\frac{1}{2}$ is of order $O \left(\theta^{\frac{1}{2}}\right)$ when $\theta \rightarrow 0$. $I_{\text{crit}, d}$ is also convergent, because $1/\zeta(s) = O(\log(\Im(s))$ by formula (3.11.8) in \cite{Titchmarsh:riemann}.\medskip

We use the following result of Valiron \cite[Theorem 9.7]{Titchmarsh:riemann}:
there exists $\alpha >0$ and a sequence $(T_k)$ such that for all $ k \in \N$, we have $ k < T_k < k + 1$  and $|\zeta(s)| > T_k^{-\alpha}$ for all $s$ 
such that  $ -1 \leq \Re(s) \leq 2$  and $|\Im(s)| = T_k$ (see the yellow rectangle in Figure \ref{fig:bande_critique}). Then, if we apply the residue theorem with the positively oriented rectangle of vertices
$\frac{3}{2} \mp i T_k$ and $-\frac{1}{2} \pm i T_k$ on the inverse Mellin transform integral, and let $k$ grow to $+\infty$, the contribution of the horizontal segments tends to $0$. As a result, the integral on the line of equation $\Re(s) = c_2$ of the Mellin transform of $\log(Zon_d (e^{-\theta} \un ))$ is equal to the sum (in the descending order regarding the real part) of the residue at 1, $I_{\text{crit}}$, the residue at 0, and the integral on the line made up of the complex numbers $s$ such that $\Re(s) = - \frac{1}{2}$. 

\begin{align*}
    \log(Zon_d (e^{-\theta} \un )) =  &\sum_{\delta=1}^{d-1}\left(p_{d,\delta}\frac{\zeta(\delta+2) \Gamma(\delta+1)}{\zeta(\delta+1) \theta^{\delta+1}} \right) + I_{\text{crit}, d} (\theta)+ C \\
    &+\frac{1}{2i\pi} \int_{ -\frac{1}{2}- i \infty}^{- \frac{1}{2} + i \infty} \frac{\mathbf{\Pi}_d (\zeta)(s) }{\zeta(s)}\zeta(s+1)\Gamma(s)  \theta^{-s} {d}s,
\end{align*}

with $C =  2 \mathbf{\Pi}_d \left(\log(2\pi) \zeta -  \zeta'\right) (0) + 2 \mathbf{\Pi}_d \left(\zeta\right)(0) \log(\theta)$.\medskip

\end{proof}



\section{PROOF OF THE THEOREM}\label{section_proof_main}

The proof of the first theorem relies on the saddle point method in $d$ dimensions (see Theorem VIII.3 in \cite{Flajolet:AC} for the one-dimensional case). Saddle point analysis consists in fetching the coefficient of the series expansion of a function by using Cauchy's integral formula. We denote $[x^n]f(x)$ the $n^{th}$ coefficient of the series expansion of $f$ with respect to $x$. Coming back to the formula (\ref{partition_function}), we extend the function $Zon_d$ on the open unit disk of $\C$ centered at 0 for each variable $x_i$. Then for any integer $i$, $1 \leq i \leq d$, given a positively oriented circle $\mathcal{C}(r) $ centered at 0, of radius $r <1$, we have: 

\begin{align}\label{Cauchy_formula}
    [x^n_i]Zon_d (\xx) = \frac{1}{2 i \pi} \int_{\mathcal{C}(r)} Zon_d (\xx) x_i^{-n-1 } dx_i.
\end{align}

The whole point of the saddle point method is to find $r$ (Lemma \ref{lem_saddle_equation_solution}) such that this Cauchy's integral is asymptotically equivalent to a gaussian integral that can be computed (Proposition \ref{definition_gittenberger}). To draw a parallel between the probabilistic approach of \cite{Bureaux:polygons} and the analytic combinatorics approach of \cite{Bodini:polyomino}, the gaussian approximation written in Proposition \ref{definition_gittenberger} is actually a local limit theorem with rate $1$ (see \cite{Bureaux:LLT} for a definition), in the probabilistic approach of \cite{Bureaux:polygons} and \cite{Bogatchev:limit}. \medskip

\subsection{Admissibility of a multivariate saddle point integral method}

The generating function is studied on the cartesian product of $d$ open unit disks contained in the complex plane, therefore for $i$ between 1 and $d$, $x_i$ is a complex number with absolute value less than 1 and we denote $x_i = e^{-\theta_i + i \alpha_i}$, with $\theta_i >0$ and $\alpha_i \in ]- \pi, \pi]$. The vector notation is then $\xx = e^{- \ttheta} e^{ i \aalpha}$. $Zon_d$ is defined on this domain, and we rewrite  the integrand of (\ref{Cauchy_formula})  as:

$$Zon_d (\xx)x_i^{-n-1 } = \frac{1}{e^{- n \theta_i}} \exp \left( \log (Zon_d (  e^{- \ttheta} e^{i \aalpha}) ) - (n+1) i \alpha_i \right). $$

In the following the notation $\xx \rightarrow \textbf{l}$ means that each component $x_i$ tends to $l_i$ and there exists a constant $\epsilon>1$, such that $\frac{x_i}{x_j} \in (1/\epsilon, \epsilon)$ for all $i \neq j$. \medskip

To ease the calculation, we introduce notations for partial derivative of $Zon_d $. The existence of such derivatives was proved in the proof of Lemma \ref{lemma_derivativ}.
We denote in the rest of the section $\aa(\xx) = \left(a_i(\xx ) \right)_{1 \leq i\leq d}$, $\BB (\xx) = \left(B_{i,j}(\xx )\right)_{1 \leq i,j \leq d}$, and $\textbf{C}(\xx) = \left(C_{i,j,k}(\xx )\right)_{1 \leq i, j, k \leq d}$ respectively for the sets of the first, second, and third order partial derivative of the logarithm of $Zon_d$ at $\xx$.
E.g. for the first order, for $1 \leq j \leq d$, we have
$$a_j(\xx) = \frac{x_j \frac{\partial  }{\partial x_j }Zon_d(\xx)}{Zon_d(\xx) }. $$
The goal of this subsection is to prove the following proposition that gives the asymptotic equivalent of the coefficient of the generating function. In analytic combinatorics terms that proposition is the \emph{H-admissibility} (see the introduction of \cite{Gittenberger:hadmi}), but we will make no use of that terminology here. \medskip

\begin{prop}\label{definition_gittenberger}
 Let $\rr = e^{- \ttheta}$ be a vector in $[0,1[^d$. Then, as $e^{- \ttheta} \rightarrow \un$, we have:
 $$[\rr^{\nn}]Zon_d(\rr) = \frac{Zon_d(e^{- \ttheta})}{e^{- \nn \cdot \ttheta} \sqrt{(2 \pi)^d \det \BB( e^{- \ttheta})}}\left( \exp\left( - \frac{1}{2}(\aa(e^{- \ttheta}) - \nn)\BB(e^{- \ttheta})^{-1}(\aa(e^{- \ttheta}) - \nn) \right) + o(1) \right), $$
 uniformly on $\N^d$.
\end{prop}\medskip

This proposition relies on a few technical lemmas about asymptotic behavior of the generating partition function, following the 5 conditions of Definition 2 in \cite{Gittenberger:hadmi} labelled by (\textbf{I}) to (\textbf{V}) therein and in the lemmas below.

\begin{lem}\label{condition1} (\textbf{I})  Let $\beta $ be a real number in the interval $ (1 + \frac{d}{3} ,1 + \frac{d}{2})$, and define the cuboid $\Delta (e^{- \ttheta} )$ in $]-\pi, \pi]^d$ as 
$$\Delta (e^{- \ttheta} ) = \left\{ \aalpha \in \R^d,\text{ for }1 \leq i\leq d, \:\: |\alpha_i| < \left(\underset{1 \leq i \leq d}{\max}( \theta_i) \right)^\beta \right\}.$$

$\BB(e^{-\ttheta})$ is positive definite, and, for all $\epsilon >1$, for all $\ttheta$ such that  $\frac{1}{\epsilon} <\frac{\theta_i}{\theta_j} < \epsilon$, we have, uniformly for $\aalpha \in \Delta (e^{- \ttheta} )$, 
$$Zon_d(e^{- \ttheta} e^{i \aalpha}) =
Zon_d(e^{- \ttheta}) 
\exp \left( i\aalpha^{\intercal} \aa(e^{- \ttheta}) -
\frac{\aalpha^{\intercal} \BB(e^{- \ttheta}) \aalpha }{2 }
\right)(1 + o(1)), \text{ as }e^{-\ttheta} \rightarrow \un.$$  
\end{lem}

\begin{proof}
Let $\ttheta \in (0,+ \infty)^d$, we compute the equivalent of the values of $\aa(e^{- \ttheta})$, $\BB (e^{- \ttheta})$, and $\textbf{C}(e^{- \ttheta})$, with Lemma~\ref{lemma_derivativ}, using Kronecker's $\delta_{ij}$ notation (and for three parameters, we write $\delta_{i,j,k} =1 $ if $i = j = k $ and 0 otherwise):

\begin{align}\label{equivalent_aa}
    \aa_i(e^{- \ttheta}) \underset{\ttheta \rightarrow \boldsymbol{0}}{\sim} \frac{2^{d-1}  \zeta(d+1)}{\zeta(d)} \: \frac{1}{\theta_i  \prod_{j = 1}^d \theta_j},
\end{align}

\begin{align}\label{equivalent_BB}
    \mathbf{B}_{i,j} (e^{- \ttheta}) \underset{\ttheta \rightarrow \boldsymbol{0}}{\sim} (1 + \delta_{i,j}) \: \frac{2^{d-1}  \zeta(d+1)}{\zeta(d)} \: \frac{1}{\theta_i \theta_j  \prod_{k = 1}^d \theta_k},
\end{align}

$$  \mathbf{C}_{i,j,k} (e^{- \ttheta}) \underset{\ttheta \rightarrow \boldsymbol{0}}{\sim} (1 + \delta_{i,j} + \delta_{j,k} + \delta_{k, i} + 2 \delta_{i,j,k})  \: \frac{2^{d-1}\zeta(d+1)}{\zeta(d)} \:\frac{1}{\theta_i \theta_j \theta_k  \prod_{l = 1}^d \theta_l}. $$

The  asymptotic equivalents (\ref{equivalent_aa}) and (\ref{equivalent_BB}) will useful later  in the following Lemma and to determine the solution of the saddle equation (Lemma \ref{lem_saddle_equation_solution}).\\
$\BB(e^{-\ttheta})$ is a symmetric matrix and the matrix of the asymptotic equivalents (\ref{equivalent_BB}) is positive definite. Therefore for $\ttheta$ small enough, $\BB(e^{-\ttheta})$ is positive definite. To conclude with the Lagrange form of Taylor's expansion theorem, we ensure that for any $t \in (0,1)$, we have   
    
$$ \mathbf{C}_{i,j,k}(e^{-\ttheta + i t \aalpha}) \alpha_i \alpha_j \alpha_k \underset{\ttheta \rightarrow \boldsymbol{0}}{\rightarrow} 0, \text{ uniformly for }\aalpha \in \Delta(\ttheta).$$
    
This is the case because $ \beta \in (1 + \frac{d}{3}, 1 + \frac{d}{2})$.
\end{proof}

\begin{lem}\label{condition2} (\textbf{II})
$|Zon_d(e^{ - \ttheta + i \aalpha})| = o \left( \frac{Zon_d(e^{-\ttheta})}{\sqrt{ \det \textbf{B}(e^{-\ttheta})}} \right)$ as $\ttheta \rightarrow \boldsymbol{0}$, holds uniformly for $\aalpha \notin  \Delta (\ttheta)$
\end{lem}

\begin{proof}
    Let $\ttheta \in (0, +\infty)^d$, and take $\aalpha \notin \Delta(\ttheta)$, and start from the following equality: \\
    
    $$\log \left(  \frac{|Zon_d(e^{ - \ttheta + i \aalpha} )|}{Zon_d(e^{ - \ttheta })} \right) = 
     - \sum_{ \vv\in \P_{d+}} 2^{d(\vv)-1} \log \left(  \frac{|1 - e^{-  \ttheta \cdot \vv + i \aalpha \cdot \vv}| }{1 - e^{-  \ttheta \cdot \vv } }\right). $$
    
    Using $ \frac{\left|1 - x e^{iy}\right|} {1 - x} = \sqrt{1+ \frac{4 x sin^2\left(y/2\right) }{(1 - x)^2}}$ for $x \in ]0,1[$, we have
    
    $$\log \left(  \frac{|Zon_d(e^{ - \ttheta + i \aalpha} )|}{Zon_d(e^{ - \ttheta })} \right) = 
     - \sum_{ \vv\in \P_{d+}} 2^{d(\vv)-2} \log \left(1 +  \frac{ 4 e^{-\ttheta\cdot \vv} \sin^2\left(\frac{\aalpha \cdot \vv}{2} \right) }{(1 - e^{-  \ttheta \cdot \vv })^2 }\right).$$
    
    We can upper bound the quotient within the logarithm in the right-hand side, whose denominator is smaller than 1, and using the fact that for $0 \leq x \leq 4$, we have $\frac{\log(5)}{4}x \leq \log(1 + x)$. We obtain:

    \begin{align}\label{formul:log_ratio_gf}
        \log \left(  \frac{|Zon_d(e^{ - \ttheta + i \aalpha} )|}{Zon_d(e^{ - \ttheta })} \right) \leq
     -  \log(5) \sum_{ \vv\in \P_{d+}} 2^{d(\vv)-2}  \left(  e^{-\ttheta \cdot \vv} \sin^2\left(\frac{\aalpha \cdot \vv}{2}\right) \right).
    \end{align}

    We denote $U_{\ttheta, \aalpha}$ as

    $$U_{\ttheta, \aalpha} =  \sum_{ \vv\in \P_{d+}} 2^{d(\vv)-2}  \left(  e^{-\ttheta \cdot \vv} \sin^2\left(\frac{\aalpha \cdot \vv}{2}\right) \right) $$
    
   To obtain the little-$o$ of the Lemma, it is sufficient to lower bound $U_{\ttheta, \aalpha}$ with a polynomial bound. Since $\aalpha \notin \Delta(\ttheta)$, there is an integer $1 \leq k \leq d$ such that $\alpha_k \geq \max_{1 \leq i \leq d }(\theta_i)$. Without loss of generality, we can suppose that $k=1$. We then consider the family of primitive vectors $\left( \textbf{v}_i = (i, 1, \textbf{0)} \right)_{i  \geq 1}$:\\
   
   $$ \sin^2(\aalpha \cdot \vv_i  ) = \sin^2(\frac{i \alpha_1 +\alpha_2 }{2}) $$
   
   We focus on the function $x \mapsto sin^2(\frac{x \alpha_1 +\alpha_2 }{2}) $ which is $\frac{ \pi}{\alpha_1}$-periodic. By construction, $\alpha_1  \leq \pi$, so we have $\frac{ \pi}{\alpha_1} \geq 1$. This final inequality guarantees that the $k$-th element of the sequence $(\vv_i)$ satisfying $sin^2(\frac{i \alpha_1 +\alpha_2 }{2} ) \geq \frac{1}{4}$ is lower that $2k$. Thus
   
   $$ U_{\ttheta, \aalpha} \geq  \sum_{i = 1}^{\infty}  e^{-\ttheta \cdot \vv_i}  \sin^2\left(\frac{\aalpha \cdot \vv}{2}\right) \geq \frac{1}{4} \sum_{i = 1}^{\infty}  e^{-2 i \underset{1 \leq j \leq d }{\max}(\theta_j) -1 }. $$
   
   Finally, elementary operations on the right term gives 
   $$ U_{\ttheta, \aalpha} \geq  \frac{e^{- \underset{1 \leq j \leq d }{\max}(\theta_j) -1}}{4(1 - e^{-2 \underset{1 \leq j \leq d }{\max}(\theta_j) })} \sim \frac{e^{-1}}{8 \underset{1 \leq j \leq d }{\max}(\theta_j) }$$
   
   which gives, for a real number $c > 0$: 
   
   $$|Zon_d(e^{ - \ttheta + i \aalpha} )| = O\left(  Zon_d(e^{ - \ttheta }) e^{- c \underset{1 \leq j \leq d }{\max}(\theta_j)^{-1}}  \right)$$
   
   The asymptotic equivalent of the determinant of $\BB(e^{-\ttheta})$ derives directly from the expression (\ref{equivalent_BB}), and is
   
   $$ \det\BB(e^{-\ttheta}) \underset{\ttheta \rightarrow \boldsymbol{0}}{\sim} \frac{2^{d-1} (d+1)  \zeta(d+1)}{\zeta(d)} \: \frac{1}{ \prod_{k = 1}^d \theta_k^3}.$$

    The big \textit{O} of the inverse of the square root of the determinant follows: 

    $$\frac{1}{\sqrt{\det \BB(e^{-\ttheta})}} = O \left(\underset{1 \leq j \leq d }{\max}(\theta_j)^{ \frac{3 d}{2}} \right),$$
    
    and we get the uniform convergence of the lemma. 
    \end{proof}

\begin{lem}\label{condtion345}
    \begin{itemize} The following properties hold:
        \item(\textbf{III}) The eigenvalues $ \lambda_1(\ttheta),..., \lambda_d(\ttheta)$ of $\textbf{B}(e^{-\ttheta})$ all tend to $+\infty$ as $\ttheta \rightarrow \boldsymbol{0}$.
        
        \item(\textbf{IV}) $\textbf{B}_{ii}(e^{-\theta}) = o (a_i(e^{-\ttheta})^2)$ as $\ttheta \rightarrow \boldsymbol{0}$.
        
        \item (\textbf{V}) For $\ttheta \rightarrow \boldsymbol{0}$, and $\aalpha \in [- \pi, \pi]^d \setminus \{\boldsymbol{0}\} $, we have  $|Zon_d(e^{-\ttheta + i \aalpha})| < Zon_d(e^{-\ttheta})$.
    \end{itemize}
\end{lem}

\begin{proof}
Property (\textbf{III}) follows from the equivalence (\ref{equivalent_BB}), Property (\textbf{IV}) directly comes from equivalences (\ref{equivalent_aa}) and (\ref{equivalent_BB}), and Property (\textbf{V}) is a direct consequence of the inequality (\ref{formul:log_ratio_gf}) that stands true for $\aalpha \in [- \pi, \pi]^d \setminus \{\boldsymbol{0}\}$. 
\end{proof}

\begin{proof}[proof of Proposition \ref{definition_gittenberger}]
With Lemmas \ref{condition1}, \ref{condition2}, \ref{condtion345}, the $d$-dimensional function $Zon_d$ satisfies all 5 conditions of Definition 2 in \cite{Gittenberger:hadmi}, therefore it satisfies Theorem 4 in \cite{Gittenberger:hadmi} which gives Proposition \ref{definition_gittenberger}.

\end{proof}

\subsection{The saddle-point equation} \label{section:saddle_equation}
 Proposition \ref{definition_gittenberger} gives an asymptotic equivalent for the coefficient of $Zon_d$. Let's consider $\nn \in (\N^*)^d$ the vector of the dimensions of the box. In order to compute $z_{d}(\nn)$ the number of  lattice zonotopes in this box, we determine $\ttheta$ as the solution of (\ref{saddle_equation}) which is often called the saddle point equation. This solution is the $\ttheta$ that cancels the exponential term in Proposition \ref{definition_gittenberger}: 

\begin{lem}\label{lem_saddle_equation_solution}
For all $\epsilon >0$, the vector $\widetilde{\ttheta}_{\nn}$ defined by 

$$ \widetilde{\theta}_{\nn, \:i} = \left(\frac{2^{d-1} \zeta(d+1)}{\zeta(d)}\right)^{\frac{1}{d+1}} \frac{\left(\prod_{j=1}^d n_i \right)^{\frac{1}{d+1}}}{n_i} $$

satisfies, as $\nn$ goes to $\infty$ such that for all $1 \leq i, j \leq d$ we have  $1/\epsilon < \frac{n_i}{n_j} < \epsilon$: 
\begin{align}\label{saddle_equation}
    \aa(e^{- \widetilde{\ttheta}_{\nn}}) =  \nn (1 + o(1) )
\end{align}
\end{lem}

\begin{proof}
Let $\ttheta \in  (0, + \infty)^d$, and consider the asymptotic equivalence coming from (\ref{equivalent_aa}): 

$$  \aa_i(e^{- \ttheta}) - n_i \underset{\ttheta \rightarrow \boldsymbol{0}}{=} \frac{2^{d-1}  \zeta(d+1)}{\zeta(d)} \: \frac{1}{\theta_i  \prod_{j = 1}^d \theta_j} -n_i + o\left( \frac{1}{\theta_i  \prod_{j = 1}^d \theta_j} \right), \text{ for $1 \leq i \leq n$.} $$

We set each $\aa_i(e^{- \ttheta}) - n_i$ to 0, and we obtain the result by computing the product of the $(\theta_j)$:

$$  \prod_{j = 1}^d \theta_j = \sqrt[d+1]{\left(\frac{2^{d-1}\zeta(d+1)}{\zeta(d)}\right)^d \frac{1}{\prod_{j =1}^{d} n_j}}.$$

Then we obtain the wanted expression by replacing the product $\prod_{j = 1}^d \theta_j$ in each equation $\aa_i(e^{- \ttheta}) - n_i = 0$.
\end{proof}

Proposition \ref{definition_gittenberger} and Lemma \ref{lem_saddle_equation_solution} imply that the number of lattice zonotopes inscribed in a box of dimensions $n \kk$ with $\kk \in (\N^{*})^d$ is 

\begin{align}\label{p_n_equivalent}
    z_{d}(  n \kk) \underset{ n \rightarrow + \infty}{\sim} \frac{Zon_d\left( e^{- \widetilde{\ttheta}_{n \kk}} \right)}{e^{- n \kk \cdot \widetilde{\ttheta}_{n \kk}} \sqrt{(2 \pi)^d \det \BB \left( e^{- \widetilde{\ttheta}_{n\kk}} \right)}} 
\end{align}

Ultimately, as Sinai did for the two-dimensional case \cite{Sinai:polygonallines}, this final asymptotic equivalence can lead to the estimate of the number of lattice zonotopes in any box $(c_1 n, ..., c_d n)$ (with positive constants $(c_i)$). Yet, a more in-depth work (analogous of what has be conducted in Section \ref{subsection_univariable_equival}) is needed to obtain an equivalent of $Zon_d\left(e^{- \widetilde{\ttheta}_{n \kk}} \right)$, so we limit our scope to the box $[0,+\infty]^d$, which leads to the parameters $\theta_1 = ... = \theta_d = \left(\frac{2^{d-1} \zeta(d+1)}{\zeta(d) n}\right)^{\frac{1}{d+1}}.$ Theorem \ref{theorem_1} follows, and we state it hereafter with detailed notations.

\subsection{Main theorem}\label{section_resultat_final}

\begin{thm}\label{theorem_1_complet}
Let $z_d(n\boldsymbol{1})$ be the number of lattice zonotopes inscribed in $[0,n]^d$. We denote $\kappa_d = \frac{2^{d-1} \zeta(d+1)}{\zeta(d)}$, and $\mathbf{\Pi}_d $ and $(p_{d,\delta})_{1 \leq \delta < d}$ respectively the operator and the coefficients defined in Definition~\ref{definition_operator}. With $\gamma$ the contour defined in Proposition \ref{prop_equivalent_univariate}, we define the polynomial $Q_d$ and the function $I_{\text{crit}, d}$ respectively by:

\begin{align}
    Q_d(X) = (d+1) \kappa_d^{\frac{1}{d+1}}  X^d +  \sum_{\delta=2}^{d-1} p_{d,\delta-1}\frac{\zeta(\delta+1) (\delta-1)!}{\zeta(\delta) }\kappa_d^{-\frac{\delta}{d+1}} X^{\delta},
\end{align}
and 

\begin{align}
    I_{\text{crit}, d}(\theta) = \frac{1}{2 i \pi } \int_{\gamma}  \frac{\mathbf{\Pi}_d (\zeta(s)) }{ \zeta(s)} \zeta(s+1) \theta^{-s} \Gamma(s) ds 
\end{align}

As $n $ grows to $+ \infty$, we have 

$$z_d(n\boldsymbol{1}) \sim  \alpha_d n ^{\beta_d} \exp \left( Q_d(n^\frac{1}{d+1}) + I_{\text{crit}, d}\left( \left(\frac{\kappa_d}{ n}\right)^{\frac{1}{d+1}} \right) \right) $$

with $\alpha_d = \frac{\kappa_d^{\frac{d}{2(d+1)} +  \frac{2}{d+1} \mathbf{\Pi}_d [\zeta](0) } \exp \left( 2\mathbf{\Pi}_d \left[\log(2\pi)\zeta -  \zeta' \right](0) \right) }{ (2\pi)^{d/2 } \sqrt{d+1}}$, and 
$\beta_d = \frac{-1}{2(d+1)}\left(d(d+2) + 4 \mathbf{\Pi}_d \left[\zeta\right] (0) \right)$.\\

Moreover, under the hypothesis that all zeros of the Riemann $\zeta$ function in the critic stripe are simple poles, $I_{crit,d}$ can be rewritten to a sum over the set $Z$ of non-trivial zeros (named hypothesis H1), that is:

\begin{align}\label{I_crit_somme}
I_{\text{crit}, d}\left( \left(\frac{\kappa_d}{ n}\right)^{\frac{1}{d+1}} \right)  = \sum_{r\in Z} Res\left(\frac{1}{\zeta(r)}\right) \mathbf{\Pi}_d [\zeta](r) \left( \frac{n}{\kappa_d} \right)^{ \frac{r}{d+1}} \zeta(r+1)\Gamma(r).  
\end{align}
\end{thm}

\begin{proof}
We rewrite (\ref{p_n_equivalent}) using the equivalent of the asymptotic equivalent of the univariate generating function of Proposition \ref{prop_equivalent_univariate}.
\end{proof}



\section{NUMERICAL CONSIDERATIONS}\label{section_numerical_discu}

In this section, we make two numerical remarks about Theorem \ref{theorem_1_complet}, in order to ease the understanding of the behavior of $z_d(n\boldsymbol{1})$ when $d$ and $n$ grow large. 

\subsection{Discussion on \texorpdfstring{$I_{\text{crit}}$}{TEXT} }
Under hypothesis H1, $I_{\text{crit},d}$ is the sum (\ref{I_crit_somme}). In fact, even without any assumption on the poles, it can still be expressed as a sum, but the terms for each zeros would be more complex. In the critical strip, the 2 first zeros of the $\zeta$ function is at $r_1 = \frac{1}{2} + i 14.1347...$ and at $r_2 = \frac{1}{2} + i 21.0220...$. Due to the exponential decrease of the function $\Gamma$ when deviating from the real line, the term $\mathbf{\Pi}_d [\zeta](r_1) \left( \frac{n}{\kappa_d} \right)^{ \frac{r_1}{d+1}} \zeta(r_1+1)\Gamma(r_1)$ is  about $10^4$ greater than the term involving $r_2$. \medskip 

Eventually, the bounds on the density of poles $r$ given by Selberg \cite{Selberg:Riemann} leads to consider that all the weight of all zeros but $r_1$ is negligible in $I_{\text{crit},d}$. We compute the approximation for the 2, 3, and 4-dimensional cases:

\begin{align}\label{numeric_I_crit}
    I_{\text{crit}, 2}\left( \left(\frac{\kappa_2}{ n}\right)^{\frac{1}{3}} \right) &\approx -1.3579 {\scriptstyle\times} &&\hspace{-0.7cm} 10^{-10}n^{1/6}\cos\left(4.7116\ln(0.6842 n)\right)\nonumber\\
    &&& - 1.4236{\scriptstyle\times}10^{-9}n^{1/6}\sin\left(4.7116\ln(0.6842n)\right),\nonumber\\
    I_{\text{crit},3}\left( \left(\frac{\kappa_3}{ n}\right)^{\frac{1}{4}} \right) &\approx -1.2325 {\scriptstyle\times}&&\hspace{-0.8cm} 
 10^{-10}n^{1/8}\cos\left(3.5337\ln(0.2777n)\right) \nonumber\\
    &&&- 1.2921 {\scriptstyle\times}10^{-9}n^{1/8}\sin\left(3.5337\ln(0.2777n)\right),\\
    I_{\text{crit}, 4}\left( \left(\frac{\kappa_4}{ n}\right)^{\frac{1}{5}} \right) &\approx -3.1764{\scriptstyle\times}&&\hspace{-0.7cm} 10^{-9}n^{1/10} \cos\left(2.8269\ln(0.1305n)\right) \nonumber\\
    &&&- 0.0628{\scriptstyle\times}10^{-9}n^{1/10}\sin(2.8269\ln(0.1305n)). \nonumber
\end{align}

To provide an order of magnitude, $I_{\text{crit}, d}\left( \left(\frac{\kappa_d}{ n}\right)^{\frac{1}{d+1}} \right)$ is smaller that $10^{-6}$  when $n < 10^{20}$ (for the first dimensions tried).\medskip

\subsection{About moving up in dimension}

For a more down-to-earth analysis of the logarithmic equivalent, i.e. the leading term of the exponential term, we can use the expansion of $\zeta(d) = 1+ 2^d + o(2^d)$ as $d $ grows large. One can see that: 

$$ \left(\frac{2^{d-1}\zeta(d+1)}{\zeta(d)} \right)^{\frac{1}{d+1}} = 2 + O\left(\frac{1}{d}\right), \:\:\:\:\:\:\:\:\text{ as } d \rightarrow + \infty.$$

Therefore, when we get that, in higher dimension, the logarithm of $z_d(n \boldsymbol{1})$ is nearly equivalent to $ 2(d+1) n^{\frac{d}{d+1}} $.\medskip

We wish to draw attention to the fact that the approach in \cite{Bureaux:polygons} and in \cite{Barany:zonotopes} focuses on polygonal lines (respectively zonotopes) beginning at $\boldsymbol{0}$ and ending at a given point whereas in our paper, we enumerate the number of lattice zonotopes in a hypercube. Given a hypercube $[-1, 1]^d$, we can split it into $2^d$ hypercubes centered at $(\pm \frac{1}{2},..., \pm\frac{1}{2})$, and view a  lattice zonotope in the hypercube as the sum of lattice zonotopes with generators in each of the square with positive first coordinate. This short explanation, also described in \cite[Theorem 6.2]{Barany:zonotopes}, explains the additional $2^{d-1}$ in the leading term of the exponential part.



\section{ESTIMATED MOMENTS OF PARAMETERS }

In this section we establish the asymptotic behavior of the first moments of parameters which can be computed with our approach. To do that, we add a variable $u$ (that we will name a parameter variable in the follow) that acts as a counting variable for the parameter. Then the partial derivative along this variable gives us the average value of the quantity under study. \\

We begin by giving a result similar to Lemma 4.3 of \cite{Barany:zonotopes} about the number of generators. That lemma gives the average number of generators of a lattice zonotope contained in a given cone and ending at a given point. In the case when the lattice zonotope is contained in a hypercube, this average can also be computed. As said before, the number of generators of a zonotope is the diameter of its graph, which gives Theorem \ref{diameter_zono}. \\

\begin{proof}[Proof of Theorem \ref{diameter_zono}]

Recall that $[ \vv ]$ denotes the vector of the absolute values of the coordinate of $\vv$. The structure of the generating function (\ref{partition_function}) is well known as it is a partition function, and each term $\left(1- e^{ - \ttheta \cdot [ \vv ] }\right)^{-1}$ yields the contribution of the generator $\vv$. When we expand it in series, the $k$-th term represents the possibility of having $k$ times the generator $\vv$:
$$\left(1- e^{ - \ttheta \cdot [ \vv ] }\right)^{-1} = \sum_{k = 0}^{+\infty}  e^{- k \ttheta \cdot [ \vv ] } $$
Therefore, we make the following modification of each term to use the variable $u$ that encodes the number of generators in a lattice zonotope:

$$ 1 + \sum_{k = 1}^{+\infty} u e^{- k\ttheta \cdot [ \vv ]}. $$

We recall the notation $e^{-\ttheta} = \left(e^{-\theta_1}, ..., e^{-\theta_d}\right)$. We call $Zon_{d, \text{gen}}$ the modified generating function defined as 

$$ Zon_{d, \text{gen}}\left(e^{- \ttheta}, u\right) =  \prod_{ \vv\in \P_{d+}} \left( 1 + \sum_{k = 1}^{+\infty} u e^{- k\ttheta \cdot  \vv} \right)^{ 2^{d(\vv)-1}}. $$

The main idea of the proof lies in the following definition of the average number of generators for lattice zonotope inscribed in $[0,n]^d$, $\mu_{gen}^n$ from the generating function (see for instance \cite[Chapter 3]{Flajolet:AC}):
$$ \mu_{gen}^n =  \frac{[\xx^{\nn}]  \left. \frac{\partial }{ \partial u}  Zon_{d, \text{gen}}(\xx, u)\right|_{u =1} }{[\xx^{\nn}] Zon_d(\xx )},$$

where $[x^n] F(x)$ (resp. $[\xx^{\nn}] F(\xx)$) denotes the coefficient of $x^n$ (resp. $\xx^\nn$) in the series expansion of $F(x)$ (resp. $F(\xx)$). All that remains is to compute the equivalent of Proposition \ref{definition_gittenberger} and Proposition \ref{saddle_equation} for  $\left.\frac{\partial }{ \partial u} Zon_{d, \text{gen}}\right|_{u =1}$. Therefore we compute the partial derivative along $u$

$$\left.\frac{\partial }{ \partial u} Zon_{d, \text{gen}}\left( e^{- \ttheta}, u\right)\right|_{u =1} 
= \left( \sum_{\vv \in  \P_{d+} }\left(  2^{d(\vv)-1} \right) e^{- \ttheta \cdot \vv } \right) Zon_{d}\left( e^{- \ttheta }\right). $$

It naturally leads to the asymptotic equivalence between the logarithm of $\left.\frac{\partial }{ \partial u} Zon_{d, \text{gen}}\right|_{u =1}$ and the one of $Zon_{d}$ when $\ttheta \rightarrow \boldsymbol{0}$, and therefore all the framework of Section \ref{section_proof_main} can be applied to $\left.\frac{\partial }{ \partial u} Zon_{d, \text{gen}}\right|_{u =1}$. Indeed, this function satisfies Lemmas \ref{condition1}, \ref{condition2}, and \ref{condtion345} as well. Hence we obtain a statement analogue to that of Proposition \ref{definition_gittenberger}, and the same parameter given by Lemma \ref{lem_saddle_equation_solution}. Denoting $\tilde{\theta}_n = \left(\frac{2^{d-1} \zeta(d+1)}{\zeta(d) n}\right)^{\frac{1}{d+1}}$, the result is 

$$\frac{[\xx^{\nn}]  \left. \frac{\partial }{ \partial u} Zon_{d,\text{gen}}(\xx , u)\right|_{u =1} }{[\xx^{\nn}] Zon_d(\xx )} \underset{n \rightarrow + \infty}{= } \left( \sum_{\vv \in  \P_{d +} }\left(  2^{d(\vv)-1} \right) e^{ - \tilde{\theta}_n \boldsymbol{1} \cdot \vv }\right)(1 + o(1)). $$

As $d(\vv)$ denotes the number of non null coordinate of $\vv$, we rewrite this sum using the Mellin inversion formula (same process as in the proof of Lemma \ref{lemma_integral_form_gg} and of Proposition \ref{prop_equivalent_univariate}), with $c> d$: 

\begin{align*}
     \sum_{\vv \in  \P_{d + } }\left(  2^{d(\vv)-1} \right) e^{ - \tilde{\theta}_n \boldsymbol{1} \cdot \vv } &=
     \frac{1}{2 i \pi }\int_{c- i\infty}^{c+ i\infty}  \sum_{\vv \in  \P_{d +} }\left(  2^{d(\vv)-1} \right)  \frac{\Gamma(s)}{\left(\tilde{\theta}_n \boldsymbol{1} \cdot \vv \right)^s} ds\\ 
     &= \frac{1}{2 i \pi }\int_{c- i\infty}^{c+ i\infty} \sum_{k = 1}^d 2^{k-1} \binom{d}{k} \left( \sum_{n \geq 1 } \frac{\binom{n-1}{k-1}}{n^s} \right) \frac{\Gamma(s)}{\zeta(s) \tilde{\theta}_n ^s} ds .
\end{align*}

With the left propagation of the integration contour used twice in this paper, we obtain

$$\sum_{\vv \in  \P_{d +} }\left(  2^{d(\vv)-1} \right) e^{ - \tilde{\theta}_n \boldsymbol{1} \cdot \vv }  \underset{ n \rightarrow + \infty }{=} \:\: \frac{2^{d-1} }{\zeta(d) \tilde{\theta}_n^{d}} + O\left(\frac{1}{\tilde{\theta}_n^{d-1}}\right) = \frac{ \sqrt[d+1]{\kappa_d}}{\zeta(d+1)} n^{\frac{d}{d+1}} + O\left(n^{\frac{d-1}{d+1}}\right),$$

with $\kappa_d =  \frac{2^{d-1} \zeta(d+1)}{\zeta(d)}$, which concludes the proof. Finally, the diameter of the graph of a zonotope is equal to its number of generators, by construction, which is to say 

$$ \mu_{\text{gen}}^n = \mu_{\text{diam}}^n.$$

\end{proof}

We can also determine another interesting property about lattice zonotopes, the estimated size of an edge of a random lattice zonotope. Depending on $\omega$ the multiplicity function of the randomly drawn zonotope defined in Section \ref{section_combi}, an edge $\ee$ is a translation of $\omega(\vv) \vv$ for a given primitive vector $\vv$. In the following proposition, we give the estimated value of $\omega(\vv)$ and its estimated variance.

\begin{prop}\label{occurences}
 The number of occurrences of a primitive generator $\vv_0$ in a lattice zonotope inscribed in $[0,n]^d$ is distributed with mean $\mu_{\text{occ}}^n$ and variance $ (\sigma^2)^n_{\text{occ}}$ such as: 
$$\mu_{\text{occ}}^n \underset{n \rightarrow + \infty}{\sim} \frac{n^{\frac{1}{d+1}}}{\kappa_d^{\frac{1}{d+1}} ||\textbf{\textit{v}}_0||_1}, \:\:\:\:\: (\sigma^2)^n_{\text{occ}}  \underset{n \rightarrow + \infty}{\sim}\left(\frac{n^{\frac{1}{d+1}}}{\kappa_d^{\frac{1}{d+1}} || \textbf{\textit{v}}_0||_1} \right)^2. $$
\end{prop}

\begin{proof}
Without loss of generality, we can choose $\vv_0 \in \P_{d+}$. As for the previous parameter, we insert a parameter variable $u$ counting for the number of occurrences of $ \vv_0$, it is substituted for:
$$ \sum_{k = 0}^{+\infty}  e^{- k  \ttheta \cdot \vv_0}  \longrightarrow \sum_{k = 0}^{+\infty} u^k e^{- k \ttheta \cdot \vv_0} .  $$

Let $Zon_{d, \text{occ}}$ be the modified generating function, that is the function that takes $\ttheta \in (0, +\infty)^d$, and returns $ Zon_d \left( e^{- \ttheta} \right)\frac{1 -  e^{- \ttheta \cdot \vv_0} }{1 - u e^{- \ttheta \cdot \vv_0} }$.

The mean and variance are asymptotically estimated like in the previous proof,  respectively

$$ \mu_{\text{occ}}^n =  \frac{[\xx^{\nn}]  \left. \frac{\partial}{\partial u} Zon_{d, \text{occ}}( \xx , u)\right|_{u =1} }{[\xx^{\nn}] Zon_d( \xx )} \:\:\:\:\text{ and }\:\:\:\:  (\sigma^2)^n_{\text{occ}} = \frac{[\xx^{\nn}]  \left. \frac{\partial}{\partial u} u \frac{\partial}{\partial u} Zon_{d, \text{occ}} (\xx , u)\right|_{u =1} }{[\xx^{\nn}] Zon_d( \xx )} - \left.\mu_{\text{occ}}^n\right.^2$$

\end{proof}

\section*{Acknowledgments}
The author wishes to very warmly thank Olivier Bodini for his ideas and remarks during all the process of this paper. Special thanks go to Philippe Marchal and Lionel Pournin for their support, their valuable comments, remarks and corrections on the first version of this paper. \\


%

\bibliographystyle{alphaurl}
\bibliography{biblio}

\end{document}